\documentclass{article}
\usepackage[b5paper, total={5.1in, 8in}]{geometry}
\usepackage[utf8]{inputenc}
\usepackage{amsmath}
\usepackage{amsfonts}
\usepackage{amsthm}
\usepackage{xcolor}
\usepackage{hyperref}
\usepackage{cleveref}
\usepackage{listings}
\usepackage{graphicx}
\usepackage{mathtools}

\usepackage{epigraph}
\usepackage{setspace}
\usepackage{comment}

\usepackage{float}
\restylefloat{table}

\usepackage{afterpage}

\graphicspath{ {images/} }

\definecolor{codegreen}{rgb}{0,0.6,0}
\definecolor{codegray}{rgb}{0.5,0.5,0.5}
\definecolor{codepurple}{rgb}{0.58,0,0.82}
\definecolor{backcolour}{rgb}{0.95,0.95,0.92}

\lstdefinestyle{mystyle}{
    backgroundcolor=\color{backcolour},   
    commentstyle=\color{codegreen},
    keywordstyle=\color{purple},
    numberstyle=\tiny\color{codegray},
    stringstyle=\color{codepurple},
    basicstyle=\ttfamily\footnotesize,
    breakatwhitespace=false,         
    breaklines=true,                 
    captionpos=b,                    
    keepspaces=true,                 
    numbers=left,                    
    numbersep=5pt,                  
    showspaces=false,                
    showstringspaces=false,
    showtabs=false,                  
    tabsize=2
}

\newtheorem{theorem}{Theorem}[section]
\newtheorem*{theorem*}{Theorem}
\newtheorem{lemma}[theorem]{Lemma}

\newtheorem*{remark}{Remark}
\theoremstyle{definition}
\newtheorem{definition}{Definition}[section]

\title{Highly Composite Numbers \\
    \doublespacing \Large MA3950 Mathematics Master Thesis \\
    Norwegian Univeristy of Science and Technology (NTNU)
    }
\author{Author: Lars Magnus Øverlier \\
        Thesis supervisor: Andrii Bondarenko}

\date{Submission date: December 15, 2022}

\begin{document}
\maketitle
\thispagestyle{empty}

\clearpage
\pagenumbering{roman}
\section*{Abstract}
\addcontentsline{toc}{section}{\protect\numberline{}Abstract}
The main result of this thesis is to show that there are only finitely many integers \(n\) such that both \(n\) and \(d(n)\) are highly composite numbers at the same time, where \(d(n)\) is the divisor function. Bertrand's postulate \cite{RamBert19} is used many times throughout the thesis and allows us to write a proof that is as simple (and as short) as possible. This thesis is meant to solve the open problem from the ``On-Line Encyclopedia of Integer Sequences" (OEIS): \href{https://oeis.org/A189394}{A189394} \cite{A189394}. 
\medbreak
The main idea for solving the problem comes from the comment in \href{https://oeis.org/A189394}{A189394};
\(n\) will contain many primes with exponent 1 when \(n\) is a large highly composite number. This implies that \(d(n)\) contains a lot of factors of \(2\). We then estimate the factor \(2^{\beta_1}\) in \(d(n)\) in terms of the largest prime in \(d(n)\) from above and from below to give us a contradiction when \(n\) is large enough. We end by finding a list of all highly composite \(n\) such that \(d(n)\) is also highly composite.

\section*{Abstrakt}
Hovedresultatet for denne oppgaven er å vise at det kun finnes endelig mange tall \(n\) slik at både \(n\) og \(d(n)\) er ``antiprimtall", hvor \(d(n)\) er divisorfunksjonen. Gjennom hele oppgaven blir Bertrands postulat \cite{RamBert19} brukt mange ganger. Dette har gjort at bevisene kan skrives så enkelt som mulig. Oppgaven skal løse det åpne problemet fra ``On-Line Encyclopedia of Integer Sequences" (OEIS): \href{https://oeis.org/A189394}{A189394} \cite{A189394}.
\medbreak
Hovedidéen for hvordan vi løser problemet kommer fra kommentaren i ~\href{https://oeis.org/A189394}{A189394};
Når \(n\) er et stort antiprimtall, vil \(n\) inneholde mange primtall med eksponent \(1\). Det vil si at \(d(n)\) inneholder mange faktorer av \(2\). Vi estimerer faktoren \(2^{\beta_1}\) i \(d(n)\) nedenfra og ovenfra i forhold til den største primtallsfaktoren i \(d(n)\) for å  få en motsigelse når \(n\) er stor nok. Vi avslutter med å finne alle antiprimtall \(n\) slik at \(d(n)\) også er et antiprimtall.
\vfill

\epigraph   {\raggedleft\large\emph{Chebyshev said it \\
                                    and I'll say it again,\\
                                    There's always a prime \\
                                    between \(n\) and \(2n\).}}
            {\large - Bertrand's postulate}

\clearpage
\section*{Acknowledgements}
\addcontentsline{toc}{section}{\protect\numberline{}Acknowledgements}
\onehalfspacing
First of all I would like to thank my supervisor, professor Andrii Bondarenko, for his valuable time, his patience with me, and our many coffee breaks during our meetings over this past year. Thank you for your persistence with me, with this thesis, and always pushing me to make it as simple as possible.

Second, my partner, my love, and my soulmate: Ella. Thank you for being there for me and your continuous support during all of our years at NTNU.

\singlespacing

\clearpage
\onehalfspacing
\tableofcontents
\singlespacing

\clearpage
\pagenumbering{arabic}
\section{Preliminaries}
This master thesis aims to prove the open problem on the OEIS \cite{A189394}: There are only finitely many positive integers \(n\) such that both \(n\) and \(d(n)\) are highly composite at the same time, where \(d(n)\) is the divisor function. In other words, for ``large enough" \(n\), assuming that both \(n\) and \(d(n)\) are highly composite we will arrive at a contradiction. We'll start by making our definitions clear.
\begin{definition}
The divisor function \(d:\mathbb{N} \to \mathbb{N} \) is defined to be the number of positive divisors of a given positive integer \(n\). It is usually represented as the sum over the divisors of \(n\)
\[
    d(n) = \sum_{d|n}1.
\]

\end{definition}

In most cases this definition does not say much on how large \(d(n)\) is in terms of the prime factorization of \(n =p_1^{\alpha_1} p_2^{\alpha_2} \dots p_r^{\alpha_r} \). Counting the number of divisors of each of the prime powers we can see that \(p_k^{\alpha_k}\) has \(1,p_k^1,p_k^2,\dots p_k^{\alpha_k}\) as its divisors: a total of \(\alpha_k + 1\) of them. Since \(d(n)\) is multiplicative we may then write 
\[
d(n) = \prod_{j=1}^r(\alpha_j + 1).
\] That is, the number of divisors of a given number does not depend on its prime factors, only on the exponents in the prime powers.

\begin{definition}
A positive number \(n\) is said to be \textit{highly composite} if the number of divisors of \(n\) is greater than the number of divisors of every number smaller than \(n\): 

\doublespacing
\centering \(n\) is \textit{highly composite} if for every \(m < n\) we have \(d(m) < d(n)\).
\end{definition}

Building on this definition, we can see that the sequence of highly composite numbers is the sequence of numbers \(n\) where \(d(n)\) obtains a new maximum. It is easy to see that if \(n=p_1^{\alpha_1} p_2^{\alpha_2} \dots p_r^{\alpha_r}\) is a highly composite number then \(n\) consists of the smallest \(r\) primes with certain exponents \(\alpha_k\). Using larger primes with the same exponents yields a larger number with the same number of divisors. Furthermore, the exponents must be in decreasing order \(\alpha_1 \geq \alpha_2 \geq \cdots \geq \alpha_r\). If not, then rearranging the exponents in a decreasing order yields a smaller number with the same number of divisors. We also have that \(\alpha_r = 1\) for all highly composite \(n\), except for the small cases when \(n=1, 4, 36\), which we will see as a consequence of the first lemma.

\section{Lemmas}
In this section we estimate the factor \(2^{\beta_1}\) in \(d(n)\) from above using \cref{upperbound2beta1}, and from below using the other lemmas. The first lemma is almost identical to Erdös' lemma 2 in \cite{Erd44}, but here it is written and proven to hold for all \(n > 50400\).

\begin{lemma} \label{altErd44L2}
Let \(n = p_1^{\alpha_1} p_2^{\alpha_2} \dots p_r^{\alpha_r}\) be a highly composite number greater than \(50400\). Then \begin{equation} \label{lemma1condition}
    \alpha_j = 1 \text{ for all } p_j \in \left[\frac{p_r+1}{2},p_r+1\right].
\end{equation}
\end{lemma}

\begin{proof}
Take any \(p_j \in [\frac{p_r+1}{2},p_r+1]\) and assume \(\alpha_j \geq 2\). Then \(\alpha_{j-2} \geq \alpha_{j-1} \geq \alpha_j \geq 2\) as mentioned in the introduction. Consider the number
\[
    n' = n \cdot \frac{p_{r+1}p_{r+2}}{p_j p_{j-1}p_{j-2}}.
\]

Clearly, \(p_{r+1}p_{r+2}\) contributes to a factor of \(4\) in \(d(n')\), and \(d(\frac{n}{p_j p_{j-1}p_{j-2}})\) is at least \(d(n) \cdot \left(\frac{2}{3}\right)^3\). So
\[
    d(n') \geq d(n)\cdot 4 \cdot \left(\frac{2}{3}\right)^3 > d(n).
\]

Our goal is to show that \(n' < n\) giving us a contradiction to \(n\) being highly composite.

Bertrand's postulate \cite{RamBert19} gives us that \(p_{k+1} < 2p_k\). So
\[
    p_{r+1}p_{r+2} < 8 p_r^2 \text{, and } \frac{1}{{p_j p_{j-1}p_{j-2}}} < \frac{8}{p_j^3}
\]
Since \(p_j \in [\frac{p_r+1}{2},p_r+1] \) we also have \(p_j \geq \frac{p_r+1}{2}\), and so
\[
    \frac{1}{p_j^3} \leq \frac{1}{(\frac{p_r + 1}{2})^3} = \frac{8}{(p_r+1)^3} < \frac{8}{p_r^3}.
\]
Thus
\[
    n' = n\frac{p_{r+1}p_{r+2}}{p_j p_{j-1}p_{j-2}} < n\frac{512}{p_r},
\]
which is less than \(n\) whenever \(p_r > 512\). This means we have a contradiction whenever the largest prime \(p_r\) in \(n\) is larger than \(512\). We check \(p_r \leq  512\) in \cref{smallcasesl21} to get that every \(n>50400\) satisfies the lemma.
\end{proof}

\begin{remark} \label{erdl2remark}
For \(n > 54000\) we get that the last two exponents \(\alpha_{r-1}\) and \(\alpha_r\) must be \(1\) by \(\pi(p_r+1)-\pi(\frac{p_r+1}{2}) \geq 2 \) for \(p_r \geq 11 = R_2\), where \(R_k\) are the so called Ramanujan primes \cite{RamBert19}. Checking a list of highly composite numbers \(\cite{Flamm03}\) when \(n \leq 50400\) gives us that \(\alpha_r = 1\) always, except for small cases when \(n = 1,4 \text{ or } 36\).\\
Also, if we write \(d(n) = 2^{\beta_1} q_2^{\beta_2} \dots q_s^{\beta_s}\) for \(n\) highly composite, we have that 

\begin{equation}
    \beta_1 \geq \pi(p_r+1)-\pi\left(\frac{p_r+1}{2}\right)
\end{equation}
where \(\pi(x) = \sum_{p \leq x} 1\) is the prime counting function. This is because each prime \(p_j \in [\frac{p_r+1}{2},p_r+1]\) of \(n\) contributes to a factor of \(2\) in \(d(n)\).
\end{remark}

The next lemma gives us the necessary upper bound for \(2^{\beta_1}\) in \(d(n)\).
 
\begin{lemma} \label{upperbound2beta1}
Suppose that \(d(n) \geq 12\) is a highly composite number and write \(d(n) = 2^{\beta_1} q_2^{\beta_2} \dots q_s^{\beta_s}\). Then
\[
    2^{\beta_1} < 8 q_s^2.
\]
\end{lemma}

\begin{remark}
Note that we don't assume \(n\) to be highly composite.
\end{remark}

\begin{proof}
We begin by choosing the smallest \(1 \leq h \leq \beta_1\) such that
\[
    2q_s 2^{\beta_1 -h} < 2^{\beta_1}.
\]
In other words \(h\) is the smallest integer such that \(2^h > 2q_s\) holds. 

If we cannot find \(h \leq \beta_1\), we set \(h = \beta_1\) so \(h-1 > \frac{\beta_1 - 1}{2}\) for \(\beta_1 \geq 2\). Then \(2^{h-1} \leq 2q_s\) and we skip the next step.

If \(h \leq \beta_1\), consider
\[
    m = d(n) q_{s+1} 2^{-h}.
\]
Then \( m < d(n) 2q_s 2^{-h}\) by Bertrand's postulate \cite{RamBert19}, which again is less than \(d(n)\) by construction of \(h\).

Then certainly we must have \(d(m) < d(d(n))\) since \(d(n)\) is highly composite. This yields
\[
    2 (\beta_1 -h + 1) < \beta_1 + 1 \text{ or } h > \frac{\beta_1+1}{2}
\]
That is, whenever \(2^h > 2q_s\) we have \(h > \frac{\beta_1+1}{2}\). Since \(h\) was chosen to be the smallest number with the property that \(2^h > 2q_s\) we then have \(2q_s \geq 2^{h-1}\) so

\begin{align*}
    2q_s    &\geq 2^{h-1} > 2^{\frac{\beta_1+1}{2}-1} \\
    4q_s^2  &> 2^{\beta_1 -1} \\
    8q_s^2 &> 2^{\beta_1}
\end{align*}
\end{proof}

The next lemma tells us that the largest prime in \(d(n)\) is bounded above by the largest prime in \(n\).

\begin{lemma} \label{qslesspr}
Write \(n\) and \(d(n)\) as earlier and suppose that both \(n\) and \(d(n)\) are highly composite. Then if \(p_r \geq 13\) we have
\begin{equation}
    q_s < p_r.
\end{equation}
\end{lemma}

Our proof is quite simple and is a just consequence of \cref{upperbound2beta1}.

\begin{proof}
We need only to prove that \(p_r > \alpha_1 + 1\) for the multiplication formula for \(d(n)\) gives us that \(\alpha_1+1 \geq q_s\). We also have that the exponents in \(n\) are decreasing \(\alpha_1 \geq \cdots \geq \alpha_r\) if \(\alpha_1 + 1\) is not prime. We now apply \cref{upperbound2beta1}, but this time on \(n\) highly composite, giving us 
\begin{equation*}
    \begin{split}
        8p_r^2 &> 2^{\alpha_1} \text{, so} \\
        p_r &> 2^{\frac{\alpha_1-3}{2}} \geq \alpha_1 + 1 \text{, if } \alpha_1  \geq 12
    \end{split}
\end{equation*}
Meaning for \(p_r \geq 13\) we have \(p_r > \alpha_1 + 1 \geq q_s\) and we are done. The result also holds for some smaller \(p_r\), but it is not needed here.
\end{proof}

The next lemma gives us a bound for \(\pi(p_r+1)-\pi\left(\frac{p_r+1}{2}\right)\) from Ramanujan's proof of Bertrand's postulate. The proof can be found in \cite{RamBert19}.

\begin{lemma}\label{ramanujanlowerbound}
\begin{equation} \label{pibound2xx}
    \pi(x) - \pi(\frac{1}{2}x) > \frac{1}{\log x}(\frac{1}{6}x - 3 \sqrt{x})\text{, for } x > 300.
\end{equation}
\end{lemma}

From these lemmas we move on to the main result.

\clearpage
\section{Main result}
\begin{theorem*}
There are only finitely many integers \(n\) such that both \(n\) and \(d(n)\) are highly composite.
\end{theorem*}

\begin{proof}
Suppose that both \(n\) and \(d(n)\) are highly composite, and write them as earlier. Since \(d(n)\) is highly composite we have the upper bound from \cref{upperbound2beta1} for \(2^{\beta_1}\)
\[
    2^{\beta_1} < 8 q_s^2.
\]

Now, assuming \(p_r \geq 2164\), the following calculation gives us a lower bound on \(2^{\beta_1}\):
\begin{equation} \label{finalcalculation}
    \begin{aligned}
        2^{\beta_1} &\geq 2^{\pi(p_r+1)-\pi\left(\frac{p_r+1}{2}\right)} &\text{(\cref{altErd44L2})}\\
                &> 2^{\frac{1}{\log (p_r+1)}(\frac{1}{6}(p_r+1) - 3 \sqrt{p_r+1})} &\text{(\cref{ramanujanlowerbound})} \\
                &> 8p_r^2 \\
                &> 8q_s^2 &\text{(\cref{qslesspr})}
    \end{aligned}
\end{equation}
Where we have used the fact that line 3 holds whenever \(p_r \geq 2164\). This gives us the contradiction.

In \cref{loweringappendix} we use a search on all primes \(p_r \leq 2164\) to find that  
\[
    2^{\pi(p_r+1)-\pi\left(\frac{p_r+1}{2}\right)} \geq 8p_r^2
\]
also holds for all primes \(181 \leq p_r < 2164\). This means there are no highly composite numbers \(n\) with largest prime \(p_r \geq 181\) such that \(d(n)\) is also highly composite.
\end{proof}

We end the main part of this thesis with \cref{listhcn}; A list of all the highly composite numbers \(n\) such that \(d(n)\) is also highly composite. This table is generated by checking if \(d(n)\) is highly composite for every highly composite number \(n\) with largest prime \(p_r \leq 181\). To guarantee that \(p_r \leq 181\), we use the bound obtained by Erdös in his proof of ``lemma 1" \cite{Erd44}. He obtains that \(n < p_r^{\pi(p_r)}\), meaning \(n < 181^{41} \leq 10^{556}\) for \(p_r \leq 181\). We do this in \cref{generatinghcn}.
\vfill
\begin{table}[H]
    \centering
    \begin{tabular}{|c|c|c|c|}
        \hline
        \(n\) & \(d(n)\) & \(n\) & \(d(n)\)\\
        \hline 
        1       & 1     & 55440             & 120 \\
        2       & 2     & 277200            & 180 \\
        6       & 4     & 720720            & 240 \\
        12      & 6     & 3603600           & 360 \\
        60      & 12    & 61261200          & 720 \\
        360     & 24    & 2205403200        & 1680 \\
        1260    & 36    & 293318625600      & 5040 \\
        2520    & 48    & 6746328388800     & 10080 \\
        5040    & 60    & 195643523275200   & 20160 \\  
        \hline
    \end{tabular}
    \caption{The cases where both \(n\) and \(d(n)\) are highly composite}
    \label{listhcn}
\end{table}
\vfill
\vfill

\clearpage
\appendix
\section{Appendix}
\renewcommand{\thesection}{A}

\subsection{Small cases in \cref{altErd44L2}} \label{smallcasesl21}
Here \(p_r \) is the largest prime in \(n\), \(p_j \) is the smallest prime greater than or equal to \(\frac{p_r+1}{2}\). The last column is the factor we multiply \(n\) by in the proof. These are the only cases when this factor is \(\geq 1\) for \(p_r \leq 512\).

\begin{table}[H]
    \centering
    \begin{tabular}{ c | c | c }
        \(p_r\) & \(p_j\) & \(\frac{p_{r+1} p_{r+2}}{p_j p_{j-1} p_{j-2}}\) \\
        \hline
         7 & 5  & 4.76667 \\
        11 & 7  & 2.10476 \\
        13 & 7  & 3.07619 \\
        17 & 11 & 1.13506 \\
        19 & 11 & 1.73247 \\
    \end{tabular}
\end{table}
We manually check the cases \(p_r \in [11,19]\) to satisfy \eqref{lemma1condition}. 
A simple check against a list \cite{Flamm03} of highly composite numbers yields 12 numbers with \(p_r = 7\) as its largest prime factor. Of these, only \(n=50400 = 2^5 \cdot 3^2 \cdot 5^2 \cdot 7\) and \(n=25200 = 2^4 \cdot 3^2 \cdot 5^2 \cdot 7\) do not satisfy \eqref{lemma1condition} and proving our lemma. For the smaller primes \(p_r = 2,3,5\) only a handful of highly composite numbers satisfies \eqref{lemma1condition}, namely \(n=2, 6, 60, 120 \text{ and } 240\), but is not needed here.
\\
Below is the Python program used for this search.

\begin{lstlisting}[language=Python]
p = list of primes < 600
table = [["p_r","p_j","p_r+1 p_r+2/p_j p_j-1 p_j-2"]]

for r in range(2,len(p)-2):
    j=0
    while p[j] < int((p[r]+1)/2):
        j+=1
    K = p[r+1]*p[r+2]/(p[j]*p[j-1]*p[j-2])
    if (K >= 1):
        table.append([p[r],p[j],K])
print(table)
\end{lstlisting}

\subsection{Lowering the largest prime in \cref{finalcalculation}} \label{loweringappendix}
Here we skip the lower bound by Ramanujan \cite{RamBert19} of \(\pi(x) -\pi(\frac{x}{2})\) to check directly when
\begin{equation}\label{skipramanujancrit}
    2^{\pi(p_r+1)-\pi\left(\frac{p_r+1}{2}\right)} \geq 8p_r^2,
\end{equation}
for all primes \(p_r\) up to \(2200\). The Python program for this search is written below and gives us that \eqref{skipramanujancrit} holds whenever \(p_r \geq 181\).

\begin{lstlisting}[language=Python]
primes = list of primes < 2200
table = [["p_r", "8p_r^2", "2^(pi(p_r+1)-pi((p_r+1)/2))"]]
def pi(x):
    return "number of primes less than x"
def pi2(x):
    return 2^(pi(x)-pi(x/2))
def f(x):
    return 8x^2
for p in primes:
    if f(p) < pi2(p+1):
        table.append(p, [f(p), pi2(p+1)])
print(table)
\end{lstlisting}

\subsection{Generating all highly composite \(n\) and \(d(n)\)} \label{generatinghcn}
We need only to check all highly composite \(n \leq 10^{556}\) as explained earlier. The 12500-th highly composite number is greater than \(10^{560}\) by direct evaluation from Flammenkamp's list  (``HCN.gz", \cite{Flamm03}). We also find that \(d(n) \leq 10^{100}\) by direct calculation of the 12500-th highly composite number \(n\). The list \texttt{hcn} in the program below is the list of all highly composite numbers less than \(10^{100}\) and we check if \(d(n)\) is in this list to generate our table.
This gives us \cref{listhcn}, which is same numbers as in \href{https://oeis.org/A189394}{A189394} \cite{A189394}!

\begin{lstlisting}[language=Python]
hcn = list of h.c.n. < 10^100
dhcn = []
def genlist(lst):
    x = [i for i in lst]
    for l in x:
        if "^" in l:
            k = l.split("^")
            i = [int(k[0])]*int(k[1])
            lst.remove(l)
            [lst.append(o) for o in i]
    lst = [int(i) for i in lst]
    return lst
f = open("HCN.gz")
lines = f.readlines()
for line in lines:
    line = line.replace("\n", "").split(" ")
    line = line[3:]
    line = genlist(line)
    dn = math.prod([i+1 for i in line])
    if dn in hcn: dhcn.append([dn,line])
print(dhcn)
\end{lstlisting}

\clearpage
\onehalfspacing
\bibliographystyle{plain}
\bibliography{main}
\addcontentsline{toc}{section}{References}

\end{document}